\documentclass[%
 aip,
 jmp,%
 amsmath,amssymb,
 reprint,%
]{revtex4-2}
\newcommand*{\conflictofinteres}[0]{\textbf{\\ \vspace{0.25cm}Conflict of interest}\textendash}
\newcommand{\References}{\textbf{\\ \vspace{0.25cm}References}}
\usepackage{hyperref}
\usepackage{graphicx}
\usepackage{dcolumn}
\usepackage{bm}
\usepackage{color}
\usepackage{parskip}
\newtheorem{definition}{Definition}
\newtheorem{corollary}{Corollary}
\newtheorem{remark}{Remark}
\begin{document}

\preprint{AIP/123-QED}

\title{On  positively divisible non-Markovian processes}

\author{Bilal Canturk}
\altaffiliation[Corresponding author: ]{bilal.cantuerk@physik.uni-freiburg.de}
\author{Heinz-Peter Breuer}%
\affiliation{ 
Institute of Physics, University of Freiburg, Hermann-Herder-Straße 3, D-79104 Freiburg, Germany
}
\affiliation{ 
EUCOR Centre for Quantum Science and Quantum Computing, University of Freiburg, Hermann-Herder-Straße 3, D-79104 Freiburg, Germany
}%

\date{\today}
\begin{abstract}
\textbf{Abstract}
There are some positively divisible non-Markovian processes whose transition matrices satisfy the Chapman-Kolmogorov 
equation. These processes should also satisfy the Kolmogorov consistency conditions, 
an essential requirement for a process to be classified as a stochastic process. Combining the Kolmogorov consistency 
conditions with the Chapman-Kolmogorov equation, we derive a necessary condition for positively divisible 
stochastic processes on a finite sample space. This necessary condition enables a systematic 
approach to the manipulation of certain Markov processes in order to obtain a positively divisible non-Markovian process. 
We illustrate this idea by an example and, in addition, analyze a classic example given by Feller in the light
of our approach.
\end{abstract}

\maketitle

\begin{quotation}
\textit{My inclination has always been to look
for general theories and to avoid computation. A
discussion I once had with Feller in a New York
subway illustrates this attitude and its limitations.
We were discussing the Markov property and I remarked that the Chapman-Kolmogorov equation
did not make a process Markovian. This statement
satisfied me, but not Feller, who liked computation
and examples as well as theory. It was characteristic of our attitudes that at first he did not believe
me but then went to the trouble of constructing a
simple example to prove my assertion.} - J. L. Doob \citep{Snell1997}
\end{quotation}

\section{\label{sec:level1} Introduction }

Markov processes are defined by two essential ingredients: $(i)$ the Chapman-Kolmogorov equation and 
$(ii)$ an initial probability distribution.\cite{Kampen2007} It may seem that these two features would determine uniquely 
whether a process is Markovian or not. However, as first noted by Doob to Feller, these two features yield
are necessary but not sufficient condition for the Markovianity of a given stochastic process.
This observation was also made by L\'{e}vy, who analyzed processes
that satisfy the aforementioned conditions but are still not Markovian.\cite{Levy1956} The first concrete example was 
proposed by Feller.\cite{Feller1959} According to Feller, these kinds of processes are rather "pathological" and, therefore, 
the conditions $(i)$ and $(ii)$  can be considered as the characteristic properties that determine Markov processes 
uniquely.\cite{Feller1972V1}  However, other examples were given later for both a countable state space 
\cite{Kampen2007} and a continuous state space. \cite{McCauley2014} 

A clear-cut understanding of the distinction between Markovian and non-Markovian processes is important 
in the theory of classical stochastic processes. \cite{Orsinger2018, Cinlar2013, Hanggi1977} 
In this context, the problem of distinguishing Markovian and non-Markovian processes is particularly relevant
for non-Markovian processes which are positively divisible 
(see Definition \ref{first_definition} below).\cite{Vacchini-Smirne-Breuer2011,Chruscinski2011a,Wissmann2015a} 
In the present paper, in contrast to McCauley's approach to the investigation of positively divisible non-Markovian 
processes which are based on stochastic differential equations \cite{McCauley2014}, we propose a necessary 
condition in terms of transition matrices for the existence of non-Markovian processes in the classical 
regime within finite state space which have the property of being positively divisible. We obtain this
condition by combining the Kolmogorov consistency conditions with the Chapman-Kolmogorov equation. 
We also argue that positively divisible non-Markovian processes are not pathological as they can be obtained by modifying 
specific Markov processes, using the necessary condition. 

The paper is organized as follows. We first recall some definitions of both 
stochastic processes and Markov processes within finite state space and introduce some conventions 
in Section \ref{section 2}. In Section \ref{section 3} we derive a necessary condition for a positively divisible 
non-Markovian process. 
In Section \ref{section 4}, we analyze a well-known example designed by Feller \cite{Feller1959},
and construct a positively divisible non-Markovian process by manipulating a certain Markov process, using our necessary 
condition. We also apply this example to a certain epidemic situation. We finally draw some conclusions
in Section \ref{section 5}.


\section {Markov processes and their characteristic properties}\label{section 2}

We assume a probability space  $\{\Omega,\mathcal{A},\mu\}$ with finite sample space $\Omega=\{\omega_{1},\omega_{2},\ldots, \omega_{n}\}$, $\sigma$-algebra $\mathcal{A}$ on $\Omega$, and probability measure $\mu:\mathcal{A}\rightarrow \mathbb{R}$. 
We further consider a stochastic process $X:\Omega \times T\rightarrow \mathbb{E}\subseteq \mathbb{R}$ for a set
of discrete time steps $T=\{t_0,t_1,\ldots\}$. For our purpose, we assume that $\mathbb{E}=\{0,1,\ldots,d-1\}$ for some positive natural number $d\in \mathbb{N}$ unless otherwise stated.  Then, the probability measure of $X_{t_{k}}$ is defined as $P\{X_{t_{k}}=j_{k}\}:=\mu (B_{j_{k}})$ such that $B_{j_{k}}=\{\omega \in \Omega \mid X_{t_{k}}(\omega)=j_{k},j_{k}\in\mathbb{E}\}$. In the following, we denote the process as  $X(T)=\{X_{t} \mid t\in T\}$. 
From a physical perspective, one of the main problems about a stochastic process is to investigate the probabilistic characterization of how the random variables $X_{t}$ at different times are related to each other. The study of this
problem is mainly based on introducing the 
joint probability distributions of finite order $m$ defined by means of 
\begin{equation}\label{joint probability of order m}
    P_{m}\{X_{t_{m}}=j_{m},X_{t_{m-1}}=j_{m-1},\ldots,X_{t_{1}}=j_{1}\}:=\mu (B_{j_{m}},B_{j_{m-1}},\ldots,B_{j_{1}}),
\end{equation}
where we have any ordered set of time instants $t_{1},t_{2},\ldots,t_{m}\in T$, satisfying 
$t_{1}\leq t_{2}\leq \ldots\leq t_{m}$, and their corresponding random variables
$X_{t_{1}}, X_{t_{2}},\ldots, X_{t_{m}} \in X(T)$.
This joint probability of order $m$ represents the probability that the time-dependent random variable $X_{t}$ takes some 
value $j_{1}$ at time instant $t_{1}$, some value $j_{2}$ at time instant $t_{2}$, etc. For each value $m=1,2,\ldots$ and 
each sequence of time instants $t_{1},t_{2},\ldots,t_{m}$, there corresponds a joint probability distribution of order $m$, and 
all together are called the family of the finite joint probability distribution of the stochastic process $X(T)$, which 
has to satisfy the following \textit{Kolmogorov consistency conditions},
\begin{enumerate}
    \item  $P_{m}\{X_{t_{m}}\in \mathbb{E},X_{t_{m-1}}\in \mathbb{E},\ldots,X_{t_{1}}\in \mathbb{E}\}=1$,
    \item $P_{m}\{X_{t_{m}}=j_{m},X_{t_{m-1}}=j_{m-1},\ldots,X_{t_{1}}=j_{1}\}\geq 0$,
    \item  $P_{m}\{X_{t_{\pi(m)}}=j_{\pi(m)},\ldots,X_{t_{\pi(1)}}=j_{\pi(1)}\}=P_{m}\{X_{t_{m}}=j_{m},\ldots,X_{t_{1}}=j_{1}\}$,
    \item  For all $1\leq k \leq m$ we have:
\begin{eqnarray}
 && \hspace{-8mm} P_{m-1}\{X_{t_{m}}=j_{m},\ldots,X_{t_{k+1}}=j_{k+1},X_{t_{k-1}}=j_{k-1},\ldots,X_{t_{1}}=j_{1}\} 
 \nonumber \\
 && \hspace{-6mm} 
 = \sum_{j_{k}\in\mathbb{E}}P_{m}\{X_{t_{m}}=j_{m},\ldots,X_{t_{k+1}}=j_{k+1},X_{t_{k}}=j_{k},X_{t_{k-1}} 
 = j_{k-1},\ldots,X_{t_{1}}=j_{1}\}.
\label{Kolmogorov consistency condition}
\end{eqnarray}
\end{enumerate}
The first two conditions are just a restatement of the conditions in the definition of a probability measure for a joint 
probability distribution of order $m$, while the third condition states that the probability distributions must be 
invariant under all permutations $\pi$ of its arguments. The fourth, and most relevant one in the context of this paper, 
means that any joint probability distribution of order $m-1$ can be obtained from the distribution of order $m$
by summation over all possible values of the random variable $X_{t_k}$ at time $t_k$.
This implies that the family of joint probability distributions of all orders must be consistent with each other in the sense that 
the marginal probability distributions of lower orders must be the same as those derived by means of the condition in 
\eqref{Kolmogorov consistency condition}. 
In the following,  we specifically refer to the condition in \eqref{Kolmogorov consistency condition} whenever we mention the term
\textit{consistency condition}. 

Recalling the conditional probabilities,
\begin{equation}
\begin{split}
     &P_{l\mid k}\{X_{t_{k+l}}=j_{k+l},\ldots,X_{t_{k+1}}=j_{k+1}\mid X_{t_{k}}=j_{k},\ldots,X_{t_{1}}=j_{1}\}\\
     &\quad \quad \quad =\frac{P_{l+k}\{X_{t_{k+l}}=j_{k+l},\ldots,X_{t_{k+1}}=j_{k+1},X_{t_{k}}=j_{k},\ldots,X_{t_{1}}=j_{1}\}}{P_{k}\{X_{t_{k}}=j_{k},\ldots,X_{t_{1}}=j_{1}\}},   
\end{split}
\end{equation}
one can express the joint probability distribution of order $m$ as
\begin{equation}\label{Joint-Conditional relation}
    \begin{split}
     & P_{m}\{X_{t_{m}}=j_{m},X_{t_{m-1}}=j_{m-1},\ldots,X_{t_{1}}=j_{1}\}\\
     &= P_{m-k\mid k}\{X_{t_{m}}=j_{m},\ldots,X_{t_{k+1}}=j_{k+1}\mid X_{t_{k}}=j_{k},\ldots,X_{t_{1}}=j_{1}\}\\
     &\quad \times P_{k}\{X_{t_{k}}=j_{k},\ldots,X_{t_{1}}=j_{1}\}.
    \end{split}
\end{equation}
\par  A stochastic process $X_{t}$ is said to be a Markov process if it rapidly forgets its past history such that its future state depends only on its present state.  This Markovian feature can be formulated in terms of conditional probabilities as 
\begin{equation}\label{Markov Condition-1}
    P_{1\mid k}\{X_{t_{k+1}}=j_{k+1}\mid X_{t_{k}}=j_{k},\ldots,X_{t_{1}}=j_{1}\}=P_{1\mid 1}\{X_{t_{k+1}}=j_{k+1}\mid X_{t_{k}}=j_{k}\},
\end{equation}
which is known as \textit{Markov condition}. It is assumed to hold for all $k=1,2,...$ and for all ordered time instants $t_{1}\leq t_{2}\leq \ldots \leq t_{k}\leq t_{k+1}$.  The Markov condition implies that the conditional probability $P_{1\mid 1}\{X_{t_{k}}=j_{k}\mid X_{t_{l}}=j_{l}\}$ for $t_{l}\leq t_{k}$ plays a central role in the theory of Markov processes. 

Before proceeding further, let us introduce some common abbreviations to simplify notation.
To this end, we define $\mathbb{R}_{+}:=\{r\in\mathbb{R}\mid r\geq 0\}$, and introduce the set of probability vectors of dimension $d$ as 
$\mathbb{P}^{d}=\{(p_{1},p_{2},\ldots,p_{d})^T\mid p_{k}\in \mathbb{R}_{+}, \sum_{i=1}^{d}p_{i}=1\}$.
Then, we define:
\begin{enumerate}
    \item [i.]  $p_{m}(j_{k},t_{k};\ldots;j_{1},t_{1}):=P_{m}\{X_{t_{k}}=j_{k},\ldots, X_{t_{1}}=j_{1}\}$ for all $t_{1}\leq \ldots \leq t_{k}\in T$ and $j_{1},\ldots,j_{k}\in \mathbb{E}$,
    \item [ii.] $p_{l\mid k}(j_{k+l},t_{k+l};\ldots;j_{k+1},t_{k+1}\mid j_{k},t_{k};\ldots;j_{1},t_{1}):=P_{l\mid k}\{X_{t_{k+l}}=j_{k+l},\ldots,X_{t_{k+1}}=j_{k+1}\mid X_{t_{k}}=j_{k},\ldots,X_{t_{1}}=j_{1}\}$,
    \item [iii.] $p(j,t \mid j',t'):=P_{1\mid 1}\{X_{t}=j \mid X_{t'}=j'\}$ so that $p(j,t \mid j',t')=p_{1\mid 1}(j,t \mid j',t')$ for all $t,t'\in T$ and for all $j,j'\in\mathbb{E}$,
    \item[iv.] $p(j,t):=P_{1}\{X_{t}=j\}$ so that $p(j,t)=p_{1}(j,t)$, 
    \item [v.] $\mathbf{p}(t):=(p(j=0,t),p(j=1,t),\ldots,p(j=d-1,t))^T\in\mathbb{P}^{d}$,
\end{enumerate}
where $\mathbf{p}(t)$ is known as one-point probability vector. The conditional probabilities $p(j,t\mid j',t')$ are also known as transition probabilities. Defining a transition matrix, $T(t,t')=(T_{jj'}(t,t'))$ with 
$T_{jj'}(t,t'):=p(j,t \mid j',t')$, one can write the evolution of the time-dependent random variable $X(T)$ in terms of the transition probabilities as 
\begin{align}
      p(j,t)&=\displaystyle \sum_{j'\in \mathbb{E}} p_{2}(j,t;j',t') \nonumber \\ 
      &=\displaystyle \sum_{j'\in \mathbb{E}} p(j,t \mid j',t')p(j',t')\nonumber \\
      &=\displaystyle \sum_{j'\in \mathbb{E}} T_{jj'}(t,t') p(j',t')\label{Evolution equation},
\end{align}
which takes on the matrix multiplication form, 
\begin{equation} \label{one-point-matrix}
 \mathbf{p}(t)=T(t,t')\mathbf{p}(t'). 
\end{equation} 
In passing, we note that a general transition matrix $T(t,t')$, by virtue of its definition, is characterized  by the following properties: 
\begin{enumerate}
    \item [1.] Its elements are non-negative, $T_{jj'}(t,t')\geq 0$,
    \item [2.] $\sum_{j}  T_{jj'}(t,t')=1$ for all $j'$,
    \item [3.] $T_{jj'}(t',t')=\delta_{j,j'}$,
    \item [4.] $p_2(j,t;j',t')=T(j,t \mid j',t')p(j',t')$ for all $t \geq t'$.
\end{enumerate}  
Note that the fourth property implies equation \eqref{one-point-matrix}.

Let us now assume a 
joint probability distribution of order $3$, $p_{3} (j_{3},t_{3};j_{2},t_{2};j_{1},t_{1})$ for any $t_{1}\leq t_{2}\leq t_{3}$ which 
corresponds to a Markov process $X(T)$. Using the Markov condition, one can write
\begin{equation}
\begin{split}
      p_{3}(j_{3},t_{3};j_{2},t_{2};j_{1},t_{1})&=p_{1\mid 2}(j_{3},t_{3}\mid j_{2},t_{2};j_{1},t_{1})p_{2}(j_{2},t_{2};j_{1},t_{1})\\
      &=p(j_{3},t_{3}\mid j_{2},t_{2})p(j_{2},t_{2}\mid j_{1},t_{1})p(j_{1},t_{1}),
\end{split}
\end{equation}
and summing over $j_{2}$ one gets
\begin{equation}
    p_{2}(j_{3},t_{3};j_{1},t_{1})=p(j_{1},t_{1})\displaystyle\sum_{j_{2}\in\mathbb{E}}p(j_{3},t_{3}\mid j_{2},t_{2})p(j_{2},t_{2}\mid j_{1},t_{1}).
\end{equation}
Dividing both sides by $p(j_{1},t_{1})$, one obtains the celebrated Chapman-Kolmogorov equation
\begin{equation}\label{Chapman_Kolmogorov equation}
    p(j_{3},t_{3}\mid j_{1},t_{1})=\displaystyle \sum_{j_{2}\in \mathbb{E}}p(j_{3},t_{3}\mid j_{2},t_{2})p(j_{2},t_{2}\mid j_{1},t_{1}).
\end{equation}
Equation (\ref{Chapman_Kolmogorov equation}) can be put into the form of the multiplication of transition matrices as
\begin{equation}\label{Chapman_Kolmogorov equation Matrix Form}
    T(t_{3},t_{1})=T(t_{3},t_{2})T(t_{2},t_{1}).
\end{equation}
We note that a complete probabilistic characterization of a Markov process is given by means of the evolution equation 
(\ref{Evolution equation}) and the Chapman-Kolmogorov equation (\ref{Chapman_Kolmogorov equation Matrix Form}).

A classical stochastic process on a finite sample space can also be formalized in a more compact form. To this end, 
denoting by $\mathbf{p}(0)$ the initial probability vector, the time evolution of the one-point distribution of the process can 
be described by a time-dependent stochastic matrix $\Lambda(t,0)$ as
 \begin{equation}\label{Stochastic evolution}
     \mathbf{p}(t)=\Lambda(t,0)\mathbf{p}(0).
 \end{equation}
 A general stochastic matrix $\Lambda(t,0)$ is characterized by two properties:
 \begin{enumerate}
     \item [1.] Its elements are non-negative, $\Lambda_{jj'}(t,0)\geq 0$,
     \item [2.] $\sum_{j} \Lambda_{jj'}(t,0) = 1$ for all $j'$.
 \end{enumerate} 
Accordingly, a stochastic matrix $\Lambda(t,0)$ can be considered as the generalization of a transition matrix by relaxing 
its third and fourth properties. We say that a stochastic matrix $\Lambda(t,0)$ is \textit{divisible} if, for any $t\geq s\geq 0$, 
one can write 
\begin{equation}\label{Stochastic evolution divisible}
    \Lambda(t,0)=\Lambda(t,s)\Lambda(s,0).
\end{equation}
If the stochastic matrix $\Lambda (s,0)$ is invertible, $\Lambda(t,s)$ can be uniquely written as 
$\Lambda(t,s)=\Lambda(t,0)\Lambda^{-1}(s,0)$. We note that $\Lambda(t,s)$ does not have to be a stochastic matrix. 
We can now formulate the definition of a P-divisible process:
\begin{definition}\label{first_definition}
A stochastic process $X(T)$  which is characterized by the divisible stochastic matrix $\Lambda(t,0)$ is called 
\textit{positively divisible (classical P-divisible)} if $\Lambda(t,0)=\Lambda(t,s)\Lambda(s,0)$ such that 
$\Lambda(t,s)$ is also a stochastic matrix for all $t\geq s\geq 0$.
\end{definition}
We stress that if $\Lambda(s,0)$ is not invertible, the stochastic matrix $\Lambda(t,s)$ that satisfies the P-divisibility
condition of this definition is not unique in general. Moreover, we see that P-divisibility reduces 
to the Chapman-Kolmogorov equation if the stochastic matrix under consideration corresponds to the transition matrix of 
the process. Therefore, we use the following convention: Whenever we speak of the \textit{transition matrix 
of a P-divisible process} we refer to a P-divisible stochastic process whose stochastic matrix 
$\Lambda(t,t')$ corresponds to its transition matrix $T(t,t')$ for any $t \geq t' \geq 0$.
  
It is evident that  Markov processes are P-divisible. As was noted at the beginning, there are some non-Markovian 
processes satisfying the Chapman-Kolmogorov equation which are P-divisible according to 
Definition \ref{first_definition}.\cite{Feller1959,McCauley2014} 
Moreover, it was shown by construction that there are also P-divisible processes which belong to the class of
semi-Markovian processes and are non-Markovian.\cite{Vacchini-Smirne-Breuer2011} In the next section, we derive a 
necessary condition in terms of transition matrices for  P-divisible non-Markovian processes. 


\section{A necessary condition for P-divisible non-Markovian processes}\label{section 3}

We consider a stochastic process $X(T)$ as defined in section \ref{section 2}, and its joint probability distribution of order 
$3$, $p_{3}(j_{3},t_{3};j_{2},t_{2};j_{1},t_{1})$, for the time steps $t_{1}\leq t_{2}\leq t_{3}$. Using the consistency condition 
and equation (\ref{Joint-Conditional relation}),  we can write
\begin{equation}
\begin{split}
        p(j_{3},t_{3}; j_{1},t_{1})&=\displaystyle\sum_{j_{2}\in\mathbb{E}}p_{3}(j_{3},t_{3};j_{2},t_{2};j_{1},t_{1})\\
        &=\displaystyle\sum_{j_{2}\in\mathbb{E}}p_{1\mid 2}(j_{3},t_{3}\mid j_{2},t_{2};j_{1},t_{1})p_{2}(j_{2},t_{2}; j_{1},t_{1})\\
        &=\displaystyle\sum_{j_{2}\in\mathbb{E}}p_{1\mid 2}(j_{3},t_{3}\mid j_{2},t_{2};j_{1},t_{1})p(j_{2},t_{2}\mid j_{1},t_{1})p(j_{1},t_{1}),
\end{split}
\end{equation} 
and dividing both sides by $p(j_{1},t_{1})$ we derive the conditional probability   
\begin{equation}\label{consisteny condition-part-1}
    p(j_{3},t_{3}\mid j_{1},t_{1})= \displaystyle\sum_{j_{2}\in\mathbb{E}}p_{1\mid 2}(j_{3},t_{3}\mid j_{2},t_{2};j_{1},t_{1})p(j_{2},t_{2}\mid j_{1},t_{1}).
\end{equation}
We note that the conditional probability  $ p(j_{3},t_{3}\mid j_{1},t_{1})$ of any stochastic process should satisfy equation (\ref{consisteny condition-part-1}), which is a consequence of the consistency condition. 
Now, we assume that the process is P-divisible, that is, the Chapman-Kolmogorov equation is satisfied. Then, we have
\begin{equation}\label{consistency-condition-part-2.2}
p(j_{3},t_{3}\mid j_{1},t_{1})
=\displaystyle\sum_{j_{2}\in \mathbb{E}}p(j_{3},t_{3}\mid j_{2},t_{2})p(j_{2},t_{2}\mid j_{1},t_{1}).
\end{equation}
Equating the right sides of equations (\ref{consisteny condition-part-1}) and (\ref{consistency-condition-part-2.2}), we find 
\begin{equation}\label{probabilistic third point equation}
    \displaystyle\sum_{j_{2}\in\mathbb{E}}p_{1\mid 2}(j_{3},t_{3}\mid j_{2},t_{2};j_{1},t_{1})p(j_{2},t_{2}\mid j_{1},t_{1}) = \displaystyle\sum_{j_{2}\in \mathbb{E}}p(j_{3},t_{3}\mid j_{2},t_{2})p(j_{2},t_{2}\mid j_{1},t_{1}),
\end{equation}
which is valid for the conditional probabilities of any P-divisible process. Introducing the following matrices,
\begin{equation}\label{transiton matrices of order three}
\begin{split}
    &Q_{1}^{(1)}=(q^{(1)}_{kl})_{d\times d}:=(p_{1\mid 2}(k,t_{3}\mid l,t_{2};0,t_{1}))_{d\times d},\\
    &Q_{1}^{(2)}=(q^{(2)}_{kl})_{d\times d}:=(p_{1\mid 2}(k,t_{3}\mid l,t_{2};1,t_{1}))_{d\times d},\\
    &\quad\vdots\\
    &Q_{1}^{(d)}=(q^{(d)}_{kl})_{d\times d}:=(p_{1\mid 2}(k,t_{3}\mid l,t_{2};d-1,t_{1}))_{d\times d},\\
    &R=(r_{kl})_{d\times d}=S_{0}=(s^{(0)}_{kl})_{d\times d}:=(p(k,t_{2}\mid l,t_{1}))_{d\times d},\\
    &S_{1}=(s^{(1)}_{kl})_{d\times d}:=(p(k,t_{3}\mid l,t_{2}))_{d\times d},\\
    &L_{1}=(l^{(1)}_{kl})_{d\times d}:=(p(k,t_{3}\mid l,t_{1}))_{d\times d},   \\
    &\mathbf{c}_{k}^{R}:=(r_{1k},r_{2k},\ldots,r_{dk})^T \quad (\mbox{$k^{th}$ column vector of R}),
\end{split}
\end{equation}
equation (\ref{probabilistic third point equation}) can be written as matrix equation:
\begin{equation}\label{Matrix multiplication equation three point}
   \begin{bmatrix}
      \quad &\mid &\quad & \mid  \\
        Q_{1}^{(1)}\mathbf{c}^{R}_{1}  &\mid &\ldots & \mid& Q_1^{(d)}\mathbf{c}^{R}_{d}\\
        \quad &\mid &\quad & \mid \\
    \end{bmatrix} = \begin{bmatrix}
        \quad &\mid &\quad & \mid  \\
        S_{1}\mathbf{c}^{R}_{1} &\mid &\ldots & \mid &S_{1}\mathbf{c}^{R}_{d}\\
        \quad &\mid &\quad & \mid  \\
    \end{bmatrix},
\end{equation}
or, equivalently,
\begin{equation}\label{Matrix multiplication equation three point second form}
    Q_{1}^{(k)}\mathbf{c}_{k}^{R}=S_{1}\mathbf{c}_{k}^{R},\quad \text{for all} \quad k=1,2,\ldots,d .
\end{equation}
The matrix multiplications $Q_{1}^{(k)}\mathbf{c}_{k}^{R}$ and $S_{1}\mathbf{c}_{k}^{R}$ are the $k^{th}$ column of the relevant matrices in equation (\ref{Matrix multiplication equation three point}). In passing, we note that the Chapman-Kolmogorov equation itself yields the equality, $L_{1}=S_{1}R$. Equation (\ref{Matrix multiplication equation three point second form}) (or equation (\ref{Matrix multiplication equation three point})) is a \textit{necessary condition} for the transition matrices of any P-divisible stochastic process $X(T)$, including Markovian and P-divisible non-Markovian processes, when the joint probability distribution of order $3$ is considered.  The matrix $Q_{1}^{(k)}$ will be called 
\textit{one-point memory transition matrix}
since it may be viewed as a transition matrix from one state to another state, conditioned on a given fixed previous state. 

A condition for the transition matrices of four time-step evolution can be obtained similar to that of three time-step evolution by following the reasoning leading to equation (\ref{probabilistic third point equation}). To this end, applying the Chapman-Kolmogorov equation to $p(j_{4},t_{4}\mid j_{1},t_{1})$ and $p(j_{3},t_{3}\mid j_{1},t_{1})$, and invoking the consistency condition, we obtain from $p_{4}(j_{4},t_{4};j_{3},t_{3};j_{2},t_{2};j_{1},t_{1})$ and $p_{2}(j_{4},t_{4};j_{1},t_{1})$ the following equation,
\begin{equation}\label{probabilistic fourth point equation}
\begin{split}
       &\displaystyle\sum_{j_{2},j_{3}\in\mathbb{E}}p_{1\mid 3}(j_{4},t_{4}\mid j_{3},t_{3};j_{2},t_{2};j_{1},t_{1})p_{1\mid 2}(j_{3},t_{3}\mid j_{2},t_{2}; j_{1},t_{1})
      p(j_{2},t_{2}\mid j_{1},t_{1})\\
       & \quad \quad = \displaystyle\sum_{j_{2},j_{3}\in \mathbb{E}}p(j_{4},t_{4}\mid j_{3},t_{3})p(j_{3},t_{3}\mid j_{2},t_{2})p(j_{2},t_{2}\mid j_{1},t_{1}).
\end{split}
\end{equation}
In addition to the definitions (\ref{transiton matrices of order three}) we introduce
\begin{equation}\label{transition matrices of order four}
\begin{split}
    &Q_{2}^{(1:1)}=(q^{(1:1)}_{kl})_{d\times d}:=(p_{1\mid 3}(k,t_{4}\mid l,t_{3};0,t_{2};0,t_{1}))_{d\times d},\\
    &\quad \vdots\\
    &Q_{2}^{(1:d)}=(q^{(1:d)}_{kl})_{d\times d}:=(p_{1\mid 3}(k,t_{4}\mid l,t_{3};0,t_{2};d-1,t_{1}))_{d\times d},\\
    &\quad\vdots\\
    &Q_{2}^{(d:d)}=(q^{(d:d)}_{kl})_{d\times d}:=(p_{1\mid 3}(k,t_{4}\mid l,t_{3};d-1,t_{2};d-1,t_{1}))_{d\times d},\\
    &S_{n-1}=(s^{(n-1)}_{kl})_{d\times d}:=(p(k,t_{n+1}\mid l,t_{n}))_{d\times d}, \quad n=1,2,3,\\
    &L_{2}=(l^{(2)}_{kl})_{d\times d}=(p(k,t_{4}\mid l,t_{1})),   \\
    &\mathbf{c}_{k}^{Q_{1}^{(n)}}:=(q_{1k}^{(n)},q_{2k}^{(n)},\ldots,q_{dk}^{(n)}) 
    \quad (\mbox{$k^{th}$ column  of $Q_{1}^{(n)}$ in equation (\ref{transiton matrices of order three})}),
\end{split}
\end{equation}
as well as the matrices
\begin{equation*}
    G_{2}^{(k)}:= \begin{bmatrix} 
        \quad &\mid &\quad & \mid  \\
        Q_{2}^{(1:k)}\mathbf{c}^{Q_{1}^{(k)}}_{1} &\mid &\ldots & \mid &Q_{2}^{(d:k)}\mathbf{c}^{Q_{1}^{(k)}}_{d}\\
        \quad &\mid &\quad & \mid  \\
    \end{bmatrix}, \quad k=1,2,\ldots,d.
\end{equation*}
Then, equation (\ref{probabilistic fourth point equation}) takes on the form
\begin{equation}\label{Matrix multiplication equation four point}
\begin{bmatrix} 
       \quad &\mid &\quad & \mid  \\
       G_{2}^{(1)}\mathbf{c}^{R}_{1} &\mid &\ldots & \mid &G_{2}^{(d)}\mathbf{c}^{R}_{d}\\
        \quad &\mid &\quad & \mid  \\
    \end{bmatrix} =  
    \begin{bmatrix} 
        \quad &\mid &\quad & \mid  \\
        S_{2}S_{1}\mathbf{c}^{R}_{1} &\mid &\ldots & \mid &S_{2}S_{1}\mathbf{c}^{R}_{d}\\
        \quad &\mid &\quad & \mid  \\
    \end{bmatrix},
\end{equation}
which  can be rewritten as
\begin{equation}\label{Matrix multiplication equation four point second form}
  G_{2}^{(k)}\mathbf{c}^{R}_{k}= \begin{bmatrix}
      \quad &\mid &\quad & \mid \\
        Q_{2}^{(1:k)}\mathbf{c}^{Q_{1}^{(k)}}_{1}  &\mid &\ldots & \mid& Q^{(d:k)}_{2}\mathbf{c}^{Q_{1}^{(k)}}_{d}\\
        \quad &\mid &\quad & \mid \\
    \end{bmatrix} \mathbf{c}_{k}^{R}=S_{2}S_{1}\mathbf{c}_{k}^{R}, \quad \mbox{for all $k=1,2,\ldots,d$}. 
\end{equation}
Equation (\ref{Matrix multiplication equation four point second form}) is a \textit{necessary condition} for the transition 
matrices of any P-divisible processes $X(T)$ when the joint probability distribution of order $4$ is relevant. In addition, we 
have the equality, $L_{2}=S_{2}S_{1}R$ $=S_{2}L_{1}$, which is the Chapman-Kolmogorov equation in the form of matrix 
multiplication. The matrix $Q_{2}^{(m:n)}$ can be called the \textit{two-point memory transition matrix} since it represents 
a transition matrix from one state to another state, conditioned on two given, fixed previous states. Using a similar
reasoning, one can obtain a necessary condition for the joint probability of order $5$ for any P-divisible 
process. To this end, we introduce the following quantities,
\begin{equation}\label{transition matrices of order five}
\begin{split}
 &Q_{3}^{(m:n:k)}=(q^{(m:n:k)}_{vl})_{d\times d}:=(p_{1\mid 4}(v,t_{5}\mid l,t_{4};s_{m},t_{3};s_{n},t_{2};s_{k},t_{1}))_{d\times d}, \\
 &\mbox{for fixed $s_{m},s_{n},s_{k}\in \mathbb{E}$},\\
    &S_{n-1}=(s^{(n-1)}_{kl})_{d\times d}:=(p(k,t_{n+1}\mid l,t_{n}))_{d\times d}, \quad n=1,2,3,4,\\
    &L_{3}=(l^{(3)}_{kl})_{d\times d}=(p(k,t_{5}\mid l,t_{1}))_{d\times d}, \\
    &\mathbf{c}_{k}^{Q_{2}^{(m:n)}}:=(q_{1k}^{(m:n)},q_{2k}^{(m:n)},\ldots,q_{dk}^{(m:n)})^T
    \quad (\mbox{$k^{th}$ column of $Q_{2}^{(m:n)}$ in equation (\ref{transition matrices of order four})}),     
\end{split}
\end{equation}
and the matrices
\begin{equation*}
    G_{3}^{(n:k)}:= \begin{bmatrix} 
        \quad &\mid &\quad & \mid  \\
        Q_{3}^{(1:n:k)}\mathbf{c}^{Q_{2}^{(n:k)}}_{1} &\mid &\ldots & \mid &Q_{3}^{(d:n:k)}\mathbf{c}^{Q_{2}^{(n:k)}}_{d}\\
        \quad &\mid &\quad & \mid  \\
    \end{bmatrix}, \quad n,k\in \{1,2,\ldots,d\}.
\end{equation*}
Then, the necessary condition for the joint probability of order $5$ in terms of the relevant transition matrices reads
\begin{equation}\label{Matrix multiplication equation of order five}
   \begin{bmatrix} 
        \quad &\mid &\quad & \mid  \\
        G_{3}^{(1:k)}\mathbf{c}^{Q_{1}^{(k)}}_{1}  &\mid &\ldots & \mid& G_{3}^{(d:k)}\mathbf{c}^{Q_{1}^{(k)}}_{d}\\
        \quad &\mid &\quad & \mid \\
    \end{bmatrix} \mathbf{c}_{k}^{R} = S_{3}S_{2}S_{1}\mathbf{c}_{k}^{R}, \mbox{for all $k=1,2,\ldots,d$}.
\end{equation}

This procedure may be extended to any desired order, generalizing correspondingly the definitions for the transition 
matrices in equations (\ref{transiton matrices of order three}), (\ref{transition matrices of order four}) and (\ref{transition 
matrices of order five}). Summarizing, we obtain our main result:
\begin{corollary}\label{Corollary-main result}
  \begin{sloppypar}
      \textit{(Necessary Consistency Condition)} Let $X(T)$ be a P-divisible stochastic process as has been considered above with finite state space $\mathbb{E}=\{0,1,\ldots,d-1\}$.  Then, the set of its transition matrices $\{S_{0},S_{1},\ldots,S_{m-2}\}$, $\{Q_{1}^{(k_{1})},k_{1}=1,2,\ldots,d\}$, $\{Q_{2}^{(k_{2}:k_{1})},k_{1},k_{2}=1,2,\ldots,d\}$,$\ldots,\{Q_{m-2}^{(k_{m-2}:\ldots:k_{1})},k_{1},\ldots,k_{m-2}=1,2,\ldots,d\}$, must satisfy the sequence of the corresponding conditions given by equations (\ref{Matrix multiplication equation three point second form}),(\ref{Matrix multiplication equation four point second form}),(\ref{Matrix multiplication equation of order five}), \ldots. Here, $m$ denotes the
number of time steps.
  \end{sloppypar}  
\end{corollary}

When the matrix $R$ is noninvertible and the process is P-divisible non-Markovian, there exist more than one matrix $S_{m-2}$ satisfying each of the equations $(\ref{Matrix multiplication equation three point second form}),(\ref{Matrix multiplication equation four point second form}),(\ref{Matrix multiplication equation of order five}), \ldots$. Therefore, it might not be possible to determine the joint probabilities uniquely. This is why we call it a necessary condition, because it is not sufficient in such cases.  Hereafter, we use the term \textit{memory transition matrices} for the matrices $Q_{m-2}^{(k_{m-2}:\ldots:k_{1})}$.  It is evident from Corollary \ref{Corollary-main result} that the Markovianity of a P-divisible process can be expressed in terms of these generalized transition matrices:

\begin{remark}\label{Remark_characterization of Markov Processes}
  \textit{(Markov Condition)} If a P-divisible stochastic process is Markovian then
  $Q_{m-2}^{(k_{m-2}:\ldots:k_{1})}=S_{m-2}$ for all $k_{m-2},\ldots, k_{1}\in\{1,\ldots,d\}$ and $m\geq3$. 
\end{remark}
 
Below, we apply the consistency criterion of Corollary \ref{Corollary-main result} to Feller's example,
\footnote{We have presented here the modified version of the example which was explored by Feller on pages $220$ and 
$423$ of his book\cite{Feller1972V1}} and use it to construct another P-divisible non-Markovian process. 


\section{Examples of P-divisible non-Markovian processes}\label{section 4}

The two P-divisible non-Markovian processes below are based on the joint probability of order $3$. Therefore, 
the necessary consistency condition of Corollary \ref{Corollary-main result} reduces to equation (\ref{Matrix multiplication equation three point second form}), and the Markov condition of Remark (\ref{Remark_characterization of Markov Processes}) leads to
\begin{equation}\label{Markov condition matrix form}
    Q_{1}^{(k)}=S_{1},\quad \mbox{for all} \quad  k=1,2,\ldots,d.
\end{equation}

\subsection{Feller's example} 

Let us consider the following sample space consisting of nine points, where six of them are the permutation 
of the numbers $1,2,3$,
\begin{equation}
    \Omega=\{(1,2,3),(1,3,2),(2,1,3),(2,3,1),(3,1,2),(3,2,1),(1,1,1),(2,2,2),(3,3,3)\}.
\end{equation}
All points are assumed to have the same probability $1/9$.
We now introduce the stochastic process $X(T)$ with three time steps, $T=\{t_{1},t_{2},t_{3}\}$, such that the random 
variable $X_{t_{k}}$ equals to the number appearing at the $k^{th}$ place of the points of the sample space. 
For example, for the point $\omega=(1,3,2)$, $X_{t_{1}}(\omega)=2, X_{t_{2}}(\omega)=3$ and $X_{t_{3}}(\omega)=1$ 
(we consider time flowing from right to left). Then, one can easily check that
\begin{equation}
    p(j_{k},t_{k})=\frac{1}{3}, \quad p_{2}(j_{k},t_{k};j_{l},t_{l})=\frac{1}{9}, \quad k,l,j_{k},j_{l}=1,2,3.
\end{equation}
It follows that while the random variables $\{X_{t_{1}},X_{t_{2}},
X_{t_{3}}\}$ are pairwise independent, they are not mutually independent, because the one-point memory conditional 
probabilities are different, such as $p_{1\mid 2}(2,t_{3} \mid 1,t_{2};1,t_{1})=0$, while  
$p_{1\mid 2}(2,t_{3} \mid 1,t_{2};3,t_{1})=1$. To proceed one step further, a triple $\{X_{t_{4}},X_{t_{5}},X_{t_{6}}\}$ 
can be defined exactly as the triple $\{X_{t_{1}},X_{t_{2}},
X_{t_{3}}\}$, but independent of it.  Continuing in this manner, one obtains the sequence of random variables 
$\{X_{t_{1}}, X_{t_{2}},\ldots, X_{t_{n}},\ldots\}$, which can be considered as the whole evolution of the stochastic process $X(T)$. Accordingly, we note that 
\begin{equation}\label{joint probability of 3k hierarchy}
\begin{split}
&p_{3k}(j_{3k},t_{3k};j_{3k-1},t_{3k-1};j_{3k-2},t_{3k-2};\ldots;j_{3},t_{3};j_{2},t_{2};j_{1},t_{1})\\
&\quad =p_{3}(j_{3k},t_{3k};j_{3k-1},t_{3k-1};j_{3k-2},t_{3k-2})\ldots p_{3}(j_{3},t_{3};j_{2},t_{2};j_{1},t_{1}), \quad  k=2,3,4,\ldots.     
\end{split}
\end{equation}
It was argued by Feller\cite{Feller1959} that this is a P-divisible non-Markovian process. Since the consecutive triples of the 
random variables of the process are independent of and similar to each other, it is sufficient to consider the evolution of the 
process for the segment $\{X_{t_{1}}, X_{t_{2}}, X_{t_{3}}\}$.  Therefore, we check if it satisfies the necessary condition 
in equation (\ref{Matrix multiplication equation three point second form}).  To this end, employing the definitions
\begin{equation*}
    \begin{split}
        &R=(r_{kl})_{3\times 3}=(p(k,t_{2}\mid l,t_{1}))_{3\times 3},\quad S_{1}= (s^{(1)}_{kl})_{3\times 3}=(p(k,t_{3}\mid l,t_{2}))_{3\times 3},\\
        &L_{1}= (l^{(1)}_{kl})_{3\times 3} =(p(k,t_{3}\mid l, t_{1}))_{3\times 3},\quad Q_{1}^{(1)}=(q^{(1)}_{kl})_{3\times 3}=(p(k,t_{3}\mid l, t_{2};1,t_{1}))_{3\times 3},\\
        &Q_{1}^{(2)}=(q_{kl}^{(2)})_{3\times 3}=(p(k,t_{3}\mid l,t_{2};2,t_{1}))_{3\times 3},\quad Q_{1}^{(3)}=(q_{kl}^{(3)})_{3\times 3}=(p(k,t_{3}\mid l,t_{2};3,t_{1}))_{3\times 3},
    \end{split}
\end{equation*}
we determine the relevant transition matrices:
\begin{equation}\label{Feller example N3}
\begin{split}
    &R=\frac{1}{3}\begin{bmatrix}
        1 & 1 &1\\
        1 & 1 &1 \\
        1 & 1 &1
    \end{bmatrix}, \; S_{1}=\frac{1}{3}\begin{bmatrix}
        1 & 1 & 1\\
        1 & 1 & 1 \\
        1 & 1 & 1
    \end{bmatrix}, \; L_{1}=\frac{1}{3}\begin{bmatrix}
        1 & 1 & 1\\
        1& 1 & 1 \\
        1 & 1 & 1
    \end{bmatrix},\\
    &Q_{1}^{(1)}=\begin{bmatrix}
        1 & 0 & 0\\
        0 & 0 & 1\\
        0 & 1 & 0
    \end{bmatrix}, \; Q_{1}^{(2)}=\begin{bmatrix}
        0 & 0 & 1\\
        0 & 1 & 0 \\
        1 & 0 & 0
    \end{bmatrix}, \; Q_{1}^{(3)}=\begin{bmatrix}
        0 & 1 & 0\\
        1 & 0 & 0 \\
         0 & 0 & 1
    \end{bmatrix}.
\end{split}
\end{equation}
First of all, it is clear that the Chapman-Kolmogorov equation is satisfied, that is, $L_{1}=S_{1}R$. The matrices $R,S_{1}$ 
and $L_{1}$ characterize the evolution for three time steps. Therefore, the process is P-divisible. Moreover, the necessary 
condition (\ref{Matrix multiplication equation three point second form}) is also satisfied,
\begin{equation}
    Q_{1}^{(k)}\mathbf{c}_{k}^{R}=\frac{1}{3}(1,1,1)= S_{1}\mathbf{c}_{k}^{R}, \quad \mbox{for all $k=1,2,3$}.
\end{equation}
On the other hand, since equation (\ref{Markov condition matrix form}) is not satisfied, that is, $Q_{1}^{(k)}\neq S_{1}$ 
for all $k=1,2,3$, the process is not Markovian but a P-divisible non-Markovian process.

The reason why Feller assessed such processes as pathological might be that, based on his own example, the initial 
probability vector $\mathbf{p}(t_{1})$ is already uniform. The process stays in this uniform probability distribution for all the 
triple time sets. In addition and more importantly, $R=S_{1}=L_{1}$. This means that the whole evolution of the process is 
essentially stationary. However, the following example can deviate from the uniform distribution and can be considered 
dynamic in the sense that it is not in a stationary state for the whole time evolution. 

\subsection{Two-state P-divisible non-Markovian process} 

We first consider a Markov process $X(T)$ with the state space $\mathbb{E}=\{0,1\}$ and the following joint probability of three discrete time steps evolution $(X_{t_{1}},X_{t_{2}},X_{t_{3}})$,
\begin{equation}\label{Two state Markov process}
\begin{split}
  p_{3}(j_{3},t_{3};j_{2},t_{2};j_{1},t_{1})&=\frac{q_{1}}{2}\delta_{j_{3},0}\delta_{j_{2},1}\delta_{j_{1},0}+\frac{q_{1}}{2}\delta_{j_{3},1}\delta_{j_{2},0}\delta_{j_{1},0}\\
  &+\frac{q_{2}}{2}\delta_{j_{3},0}\delta_{j_{2},1}\delta_{j_{1},1}+\frac{q_{2}}{2}\delta_{j_{3},1}\delta_{j_{2},0}\delta_{j_{1},1} , \quad j_{3}, j_{2},j_{1}\in\mathbb{E}, 
\end{split}
\end{equation}
where $\delta_{k,l}$ is the Kronecker delta function and $\mathbf{p}(t_{1})=(q_{1},q_{2})^T\in\mathbb{P}^{2}$ is the initial probability vector. Like in Feller's example, we assume that the process is a sequence of the triple $(X_{t_{1}},X_{t_{2}},X_{t_{3}})$ such that each $\mathbf{p}(t_{3k-2})$ is some arbitrary probability vector for $k=1,2,3,\ldots$. We also note that equation (\ref{joint probability of 3k hierarchy}) is valid for this example. In the following, we consider $\mathbf{p}(t_{3k-2})=\mathbf{p}(t_{1})$ for all 
$k=1,2,3,\ldots$ for simplicity. In order to characterize the process, it is sufficient to analyze the first three time steps 
of the process. One can directly check that the transition matrices are
\begin{equation}
    \begin{split}
        &R=(r_{kl})_{2\times 2}=(p(k ,t_{2}\mid l,t_{1})_{2\times 2}=\frac{1}{2}
        \begin{bmatrix}
            1 & 1 \\
            1 & 1
        \end{bmatrix},  \\
        &S_{1}=(s^{(1)}_{kl})_{2\times 2}=(p(k,t_{3}\mid l,t_{2}))_{2\times 2}=\begin{bmatrix}
            0 & 1\\
            1 & 0
        \end{bmatrix},\\
        &L_{1}=(l^{(1)}_{kl})_{2\times 2}=(p(k,t_{3}\mid l,t_{1}))_{2\times 2}=\frac{1}{2}
        \begin{bmatrix}
            1 & 1 \\
            1 & 1
        \end{bmatrix}, \\\
        & Q_{1}^{(1)}=Q_{1}^{(2)}=S_{1}.
    \end{split}
\end{equation}
The process is Markovian since equation (\ref{Markov condition matrix form}) is satisfied, i.e., $Q_{1}^{(1)}=Q_{1}^{(2)}=S_{1}$. For $\mathbf{p}(t_{1})=(0.25,0.75)^T$, we estimate the time average and the variance of the process, 
respectively, as
\begin{eqnarray}
    \mu &=& \frac{1}{3}\langle X_{t_{3}}+X_{t_{2}}+X_{t_{1}}\rangle \nonumber \\
    &=& \frac{1}{3}\sum_{j_{1},j_{2},j_{3}\in \mathbb{E}}(j_{3}+j_{2}+j_{1})p_{3}(j_{3},t_{3};j_{2},t_{2};j_{1},t_{1}) \nonumber \\
    &=& \frac{1}{3}(2q_{2}+q_{1})\approx 0.583, \nonumber \\
    \mathrm{Var}(X) &=& \left\langle\left(\frac{1}{3}(X_{t_{3}}+X_{t_{2}}+X_{t_{1}})-\mu\right)^{2}\right\rangle \nonumber \\
    &=& \sum_{j_{1},j_{2},j_{3}\in \mathbb{E}}
    \left(\frac{1}{3}(j_{3}+j_{2}+j_{1})-\mu\right)^{2}p_{3}(j_{3},t_{3};j_{2},t_{2};j_{1},t_{1})
    \nonumber \\
    &=& \frac{1}{9}(q_{1}-q_{1}^{2})\approx 0.021.
\end{eqnarray}

Now, we modify the Markov process to obtain a P-divisible non-Markovian process $\tilde{X}(T)$. To this end, we make use of the necessary condition in equation (\ref{Markov condition matrix form}), which states that a P-divisible process is not Markovian if its one-point memory transition matrices are different from each other. Accordingly, we modify the one-point memory transition matrix $Q_{1}^{(1)}$ as 
\begin{equation}\label{modification of Markov process}
    \tilde{Q}_{1}^{(1)}=\varepsilon I+(1-\varepsilon) Q_{1}^{(2)}
\end{equation}
 with $\varepsilon \in [0,1]$. The new process will not be Markovian since the Markov condition in equation (\ref{Markov condition matrix form}) is not satisfied anymore unless $\varepsilon = 0$.   $\tilde{Q}_{1}^{(1)}$ implies that the future of the process will be affected by the past if the initial value of the process is zero. 

Now, using the new $\tilde{Q}_{1}^{(1)}$ and the matrices $\{R,L_{1},Q_{1}^{(2)}\}$, we construct the following joint probability distribution of three-time-step evolution of the stochastic process $\tilde{X}(T)$,
\begin{equation}\label{Two state non-Markovian process}
    \begin{split}
         \tilde{p}_{3}(j_{3},t_{3};j_{2},t_{2};j_{1},t_{1})&=\frac{\varepsilon}{2}q_{1}\delta_{j_{3},0}\delta_{j_{2},0}\delta_{j_{1},0}+\frac{1-\varepsilon}{2}q_{1}\delta_{j_{3},0}\delta_{j_{2},1}\delta_{j_{1},0}\\
         &+\frac{1-\varepsilon}{2}q_{1}\delta_{j_{3},1}\delta_{j_{2},0}\delta_{j_{1},0}+\frac{\varepsilon}{2}q_{1}\delta_{j_{3},1}\delta_{j_{2},1}\delta_{j_{1},0}\\
         &+\frac{q_{2}}{2}\delta_{j_{3},0}\delta_{j_{2},1}\delta_{j_{1},1}+\frac{q_{2}}{2}\delta_{j_{3},1}\delta_{j_{2},0}\delta_{j_{1},1} , \quad j_{3}, j_{2},j_{1}\in\mathbb{E}.
    \end{split}
\end{equation}
\begin{figure}[h!]
    \centering
    \includegraphics[width=0.75\textwidth,height=0.5\linewidth]{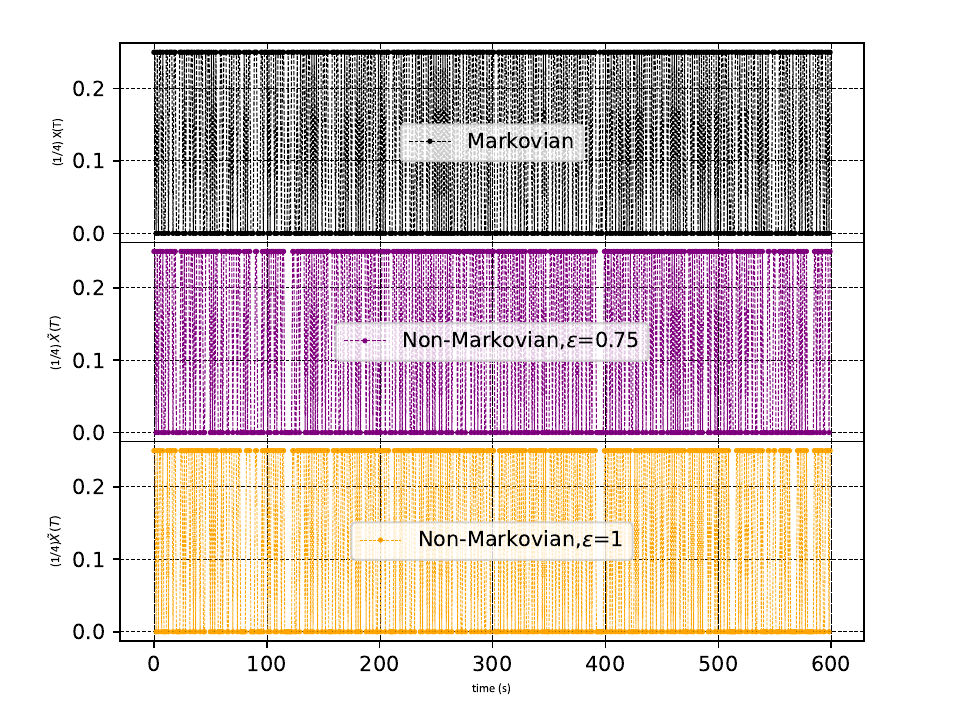}
    \caption{Realizations of the Markov process $X(T)$ with the joint probability distribution in equation (\ref{Two state Markov process}), and the non-Markovian process $\tilde{X}(T)$ with the joint probability distribution in equation (\ref{Two state non-Markovian process}) for different values of the parameter $\varepsilon$. The probability vectors $\mathbf{p}(t_{3k-2})$ have been chosen as $(0.25,0.75)^T$ for all $k=1,2,3,\ldots$. The first two time steps of the realizations have been generated using the same random numbers since there is no difference between the Markovian and the 
    non-Markovian process for the first two time steps. From these realizations, one finds the following numerical values: 
    $\mu = 0.575, \tilde{\mu}=0.585$ (for $\varepsilon =0.75$), $\tilde{\mu}=0.580$ (for $\varepsilon = 1$), $\rm{Var}(X)=0.022, \rm{Var}(\tilde{X})=0.039$ (for $\varepsilon =0.75$), $\rm{Var}(\tilde{X})=0.050$ (for $\varepsilon =1$), which 
    are in accordance with the analytical predictions within statistical error.}
    \label{figure 2}
\end{figure}
From equation (\ref{Two state non-Markovian process}), one can straightforwardly derive the modified transition matrix $\tilde{S}_{1}$ as
\begin{equation}
\tilde{S}_{1}= 
    \begin{bmatrix}
    \varepsilon q_{1} & (1-\varepsilon) q_{1}+q_{2}\\
    (1-\varepsilon) q_{1}+q_{2} & \varepsilon q_{1}
    \end{bmatrix}=q_{1}(\varepsilon I +(1-\varepsilon)S_{1})+q_{2}S_{1}.
\end{equation}
It can be checked that the new process with the transition matrices 
$\{\tilde{Q}_{1}^{(1)}, Q_{1}^{(2)}, \tilde{S}_{1}, R\}$ satisfies the necessary condition in equation (\ref{Matrix multiplication 
equation three point second form}). 
Moreover, the Chapman-Kolmogorov equation $L_1=\tilde{S}_1R$ is fulfilled, while the Markov condition 
in equation (\ref{Markov condition matrix form}) is not. 
Therefore, the new process $\tilde{X}(T)$ is a P-divisible non-Markovian process.  
Similar to Feller's example, while the random variables $\{\tilde{X}_{t_{1}},\tilde{X}_{t_{2}},
\tilde{X}_{t_{3}}\}$ are pairwise independent, they are not mutually independent. We also note that $\tilde{S}_{1}$ depends on the initial probability as was argued by H\"anggi and Thomas \cite{Hanggi1977} and by McCauley.\cite{McCauley2014}

To get an idea of how the resulting P-divisible non-Markovian process differs from the original Markov process, we show some realizations in Figure \ref{figure 2} for different values of the parameter $\varepsilon$.  The average value $\tilde{\mu} $ of $\tilde{X}(t)$ for all values of $\varepsilon$ is surprisingly equal to that of $X(t)$. However, the variance $Var(\tilde{X})$  of $\tilde{X}(t)$  is equal to $(1/9)((1+\varepsilon)q_{1}-q_{1}^{2})$, which differs from that of $X(t)$. For the initial probability, $\mathbf{p}(t_{1})=(0.25,0.75)^T$,  it is equal to $0.041$ and $0.049$ for $\varepsilon = 0.75$ and $1$ respectively.  The difference between the Markovian and the non-Markovian processes in terms of variance can also be captured partially from the realizations in Figure \ref{figure 2} by observing the gaps around, for example, the time step $t=400$. The gaps can be interpreted as an impact of memory on future states.     

As $\varepsilon \rightarrow 0$, the process approaches the original Markov process. On the other hand, the memory of the process increases by increasing $\varepsilon$, leading to an increase of the degree of non-Markovianity. In addition, it is noteworthy that the process reduces to the P-divisible non-Markovian process given by van Kampen\footnote{See the sixth example on p.79 of van Kampen's book\cite{Kampen2007}} in the limit $\varepsilon \rightarrow 1$ and for $\mathbf{p}(t_{1})=(0.5,0.5)^T$. Thus, we conclude that Markovian and P-divisible non-Markovian processes are essentially different from each other. Furthermore,  considering also equation (\ref{modification of Markov process}), this example suggests that a P-divisible non-Markovian process might arise from disturbing a certain Markov process. This disturbance induces memory effects in the process, but the memory only lasts for a short time and the process returns to its original Markovian evolution.  The memory effect manifests itself in $\mathrm{Var}(\tilde{X})$ being greater than $\mathrm{Var}(X)$. 

The two-state P-divisible non-Markovian process presented above might seem purely a mathematical artefact, but it has 
the potential to be applied to some physical problems. Below, we will present such an application.

\subsection{Physical model of the two-state P-divisible non-Markovian process} 

We consider a city with a population of $N$. Suppose a new epidemic arises in the city, and in compliance with regulations, everyone must wear a mask until a vaccine is developed. Fortunately, a vaccine is developed, and the city government plans a vaccination program with specific appointment times for the three doses. Then, the government announces that individuals who receive at least two vaccinations may remove their masks and resume normal activities and that those who have not yet received two vaccinations must continue to wear masks until the next vaccination schedule. 

\par In such a model, depending on the first dose, some individuals might not take the third dose because two doses are enough to resume their normal life. Therefore, the memory of the first dose will affect receiving the third vaccine. In addition, there might be some people who are against vaccination, so they may resist getting vaccinated. The individual's state of being vaccinated at each of three appointments can be represented with the state space $\mathbb{E}=\{0,1\}$, such that $0$ and $1$ correspond to the events "not getting vaccinated" and "getting vaccinated", respectively. We assume that each individual is vaccinated independently of each other. Then, the vaccination process of each individual at the \textit{first schedule} can be represented with the first triple steps, $\{\tilde{X}_{t_{1}},\tilde{X}_{t_{2}},\tilde{X}_{t_{3}}\}$, of the two-state P-divisible non-Markovian process given above conditioned on the following specific statements related to the memory of the first vaccination:
\begin{enumerate}
    \item [1.] Each individual gets the first vaccination with the initial probability, $\mathbf{p}(t_{1}) =(q_{1},q_{2})^T$. 
    \item [2.] Depending on vaccination status at one of the appointments, people are neutral against one of the other vaccines because they might think that they do not have to get vaccinated three times to be free of masks, and therefore, they might skip one of the vaccination appointments. We quantify this situation with the two-point conditional probabilities: $p(j,t_{3}\mid k,t_{2})=p(j,t_{3}\mid k,t_{1})=p(j,t_{2}\mid k, t_{1})=1/2$.  
    \item [3.] If an individual gets the first two vaccinations, they do not get the third one because they are free of masks: $p_{1\mid 2}(0,t_{3}\mid 1,t_{2};1,t_{1})=1$.
    \item [4.] If an individual gets the first but not the second vaccination, they get the third one surely to be free of mask: $p_{1\mid 2}(1,t_{3}\mid 0,t_{2};1,t_{1})=1$.
    \item [5.] If an individual does not get vaccinated at the first appointment, they do not also get vaccinated at the third appointment conditioned on the second vaccination with the following probabilities: $p_{1\mid 2}(0,t_{3}\mid 0,t_{2},0,t_{1})=\varepsilon$, $p_{1\mid 2}(0,t_{3}\mid 1,t_{2},0,t_{1})=1-\varepsilon$. Such an individual might be against vaccination at the beginning and later might change their opinion because of some reason, such as being informed about the benefits of vaccination. 
\end{enumerate}
The joint probability of the three vaccinations for each individual in the city whose individuals show the above characteristics can be reasonably represented by equation (\ref{Two state non-Markovian process}).  Now, our question is as follows: 

\textit{Question:} How many people will be free of masks after the first three vaccinations? 

According to the government announcement, those who have been vaccinated at least two times will be free of masks.  Therefore, those individuals will be free of masks whose vaccination state matches with one of the following events: $(1,t_{3};1,t_{2};1,t_{1})$, $(0,t_{3};1,t_{2};1,t_{1})$, $(1,t_{3};0,t_{2};1,t_{1})$ and $(1,t_{3};1,t_{2};0,t_{1})$. Then, the number of individuals who are free of masks is
\begin{equation}\label{Vaccinated people in the first schedule}
\begin{split}
        &N\big(p_{3}(1,t_{3};1,t_{2};1,t_{1})+p_{3}(0,t_{3};1,t_{2};1,t_{1})+p_{3}(1,t_{3};0,t_{2};1,t_{1})+p_{3}(1,t_{3};1,t_{2};0,t_{1})\big)\\
        &= N\left(0+\frac{q_{2}}{2}+\frac{q_{2}}{2}+\frac{\varepsilon}{2}q_{1}\right)
        =N\left(q_{2}+\frac{\varepsilon}{2}q_{1}\right)
        =N(1-q_{1}(1-\varepsilon/2)).
\end{split}
\end{equation}
In this model, it seems that the memory parameter $\varepsilon$  refers to the situation of raising awareness in those individuals who have not been vaccinated at the first appointment. These individuals might be either those against vaccination or those who have decided not to take the vaccination according to the information given in items $2$ and $5$. 

The number of unvaccinated individuals for the next schedule will be $N-N(1-q_{1}(1-\varepsilon/2))= Nq_{1}(1-\varepsilon/2)$. The vaccination process of these unvaccinated individuals for the next schedule can be represented by the second triple steps, $\{\tilde{X}_{t_{4}},\tilde{X}_{t_{5}},\tilde{X}_{t_{6}}\}$, of the two-state P-divisible non-Markovian process under the same conditions given above. However, the initial probability might be updated as $\mathbf{p}(t_{4})=(q_{1},q_{2})^T\rightarrow (q_{1}(1-\varepsilon/2),q_{2}+q_{1}\varepsilon/2)^T$ to take into account the effect of the memory on the unvaccinated individuals during the first vaccination process. Such a modification can be made if the government has recorded the vaccinated people in equation (\ref{Vaccinated people in the first schedule}) during the first schedule and used this information to update the initial probability for the second schedule. We point out that this modification would change the mean value and the variance of the P-divisible non-Markovian processes because
the equation $\mathbf{p}(t_{3k-2}) = \mathbf{p}(t_{1})$ no longer holds for all $k=1,2,3,\ldots$.  

\par We could also model the epidemic according to the Markov process whose joint probability is given by equation (\ref{Two state Markov process}). In that case, there would not be any memory effect, $\varepsilon =0$.  The number of individuals who were free of masks after the first schedule would be $N(1-q_{1})$, which is less than that of the non-Markovian case.  The number of unvaccinated individuals for the next schedule would be $Nq_{1}$, which is greater than that of the non-Markovian case. A possible lesson we can take away from this example is that in bad cases like epidemics,  raising people's awareness affects the evolution of the case in time.   

\section{Conclusion}\label{section 5}
P-divisibility is an essential property for the differential analysis of classical stochastic processes. For instance, the 
Kolmogorov forward and backward equations require the P-divisibility of the process.\cite{Feller1972V1, Feller1970V2}   
As has been stated in Corollary \ref{Corollary-main result}, we have established a necessary condition for the transition 
matrices of P-divisible stochastic processes within a finite sample space. It has some useful features and advantages as
it provides a systematic way of constructing a P-divisible non-Markovian process from a given Markovian 
process. This approach is based on the (multipoint) \textit{memory transition matrices} by which one can encode certain 
memory effects, as has been expressed by equation (\ref{modification of Markov process}).    

We have illustrated Corollary \ref{Corollary-main result} by two examples in section \ref{section 4}. First, 
Feller's example \cite{Feller1959} has been shown to be a consistent P-divisible non-Markovian process.  Therefore, in 
contrast to Feller, we have argued that P-divisible non-Markovian processes are not pathological, as they seem
to arise by injecting short-term memory into certain Markov processes.  
In the second example and its application to the epidemic, we have illustrated this idea through our construction of the 
two-state P-divisible non-Markovian process utilizing the necessary condition in equation (\ref{Matrix 
multiplication equation three point second form}). The process is distinguished by a parameter $\varepsilon$, which 
characterizes the degree of non-Markovianity and can be measured in terms of the variance of the process. 

Concluding, we have clarified that there are P-divisible non-Markovian processes which can be systematically 
constructed from certain Markov processes by taking into account so-called memory transition matrices
describing the impact of memory effects.

\begin{acknowledgments}
B.C.  gratefully acknowledges support from the Georg H. Endress foundation.
B.C. also thanks Emre G\"ursoy for helpful discussions on the manuscript.
\end{acknowledgments}

\conflictofinteres The authors declare that they have no conflicts of interest.
\References
\bibliography{Reference}

\end{document}